\documentclass[3p,preprint]{elsarticle}
\usepackage{amsthm}
\usepackage{graphicx}
\usepackage{amssymb}
\usepackage{url}
\usepackage{algorithm}
\usepackage{algpseudocode}
\usepackage{amsmath}
\usepackage{bbm}
\usepackage{graphicx}
\usepackage{amssymb}
\usepackage{url}

\usepackage{algorithm}
\usepackage{algpseudocode}

\usepackage{amsmath}

\usepackage{bm}

\DeclareMathOperator*{\argmax}{argmax}
\usepackage{color}

\newcommand{\mK}{\mathsf{K}}

\usepackage{multirow}

\usepackage{color}

\journal{Journal of Computational and Applied Mathematics}

\begin{document}

\begin{frontmatter}

\title{Parameter Tuning in the Radial Kernel-Based Partition of Unity Method by Bayesian Optimization}

\author[address-TO,address-GNCS]{Roberto Cavoretto} 
\ead{roberto.cavoretto@unito.it}

\author[address-TO,address-GNCS]{Alessandra De Rossi}
\ead{alessandra.derossi@unito.it}

\author[address-TO,address-GNCS]{Sandro Lancellotti\corref{corrauthor}}
\ead{sandro.lancellotti@unito.it}

\author[address-TO]{Federico Romaniello}
\ead{federico.romaniello@unito.it}

\cortext[corrauthor]{Corresponding author}

\address[address-TO]{Department of Mathematics \lq\lq Giuseppe Peano\rq\rq, University of Torino, via Carlo Alberto 10, 10123 Torino, Italy}

\address[address-GNCS]{Member of the INdAM Research group GNCS}

\begin{abstract}
In this paper, we employ Bayesian optimization to concurrently explore the optimal values for both the shape parameter and the radius in the partition of unity interpolation using radial basis functions. Bayesian optimization is a probabilistic, iterative approach that models the error function through a progressively self-updated Gaussian process. Meanwhile, the partition of unity approach harnesses a meshfree method, allowing us to significantly reduce computational expenses, particularly when considering a substantial number of scattered data points. This reduction in computational cost is achieved by decomposing the entire domain into several smaller subdomains, each of them with a variable radius. We provide an estimation of the complexity of our algorithm and carry out numerical experiments to illustrate the effectiveness of our approach, dealing with test and real-world datasets. 

\end{abstract}


\end{frontmatter}
\section{Introduction}

Over the past few decades, radial basis function (RBF) approximation and interpolation have emerged as a dynamic and significant tool for advancing meshfree techniques in the solution of various types of scientific and engineering problems, see e.g. \cite{bia18,Fasshauer,fra18}. They offer several advantageous features, such as straightforward implementation in higher dimensions, adaptability to various geometric configurations, and reasonable convergence properties, to name a few \cite{wen05}. However, they may lead to a full, computationally expensive and ill-conditioned linear system. To overcome this drawback, in this work we focus on a meshfree method, known as the RBF partition of unity method (RBF-PUM), which makes use of local RBF approximants accumulating all the local contributions in a global partition of unity fit. A first version of PUM is introduced in \cite{shep68} to reconstruct a function from scattered data points. This approach hinges on the concept of localizing the approximation process by a decomposition of the original big problem into several small subproblems, thus finding application in many fields of computational mathematics and scientific computing \cite{bab97,cav21a,Sandro}.  Indeed, the PUM is used to efficiently split the data within smaller subdomains or balls. 
The first combination of the PUM with the RBF interpolation goes back to \cite{wen02}, where an error analysis is also given for functions in the native space of the underlying RBFs. The RBF-PUM  proposed in this paper is obtained by a weighted sum of local RBF interpolants depending on a shape (or scale) parameter $\varepsilon$ and a variable radius $\delta$ in each subdomain.



Moreover, in this study we employ a well-studied statistical method known as {\sl Bayesian Optimization} (BO)  \cite{Practical} to simultaneously search for the optimal values of $(\varepsilon, \delta)$ within each subdomain of RBF-PUM. The BO, originally developed in the field of machine learning for optimizing complex or hard-to-assess functions, finds utility in hyperparameter tuning tasks by circumventing the need to compute and evaluate the approximations for parameter combinations that are far from optimal. All the algorithms involved in the main procedure are described and analysed in detail in order to show how this approach leads to a substantial reduction in terms of computational time. Numerical experiments on some benchmark test cases and real-world datasets such as Tonga Trench and Franke's glacier ones point out that also in applied contexts a good accuracy of the interpolant is preserved.


The paper is organised as follows. In Section \ref{rbfpum}, the RBF-PUM interpolation problem is stated. In Section \ref{bo}, BO Gaussian processes and acquisition functions are presented. Section \ref{algorithmsandcc} contains a description of the algorithms and their complexity analysis. In Section \ref{exp} numerical experiments show the efficacy of our scheme by solving interpolation problems on some test examples and real-world applications.


\section{RBF-PUM Interpolation}\label{rbfpum} 

In this section, we introduce the interpolation problem and the basic theory on RBF-PUM, highlighting the reasons that inspired this paper. 

\subsection{The RBF Method}

Let $ X = \{  \boldsymbol{x}_i, i = 1,  \ldots , N \}$ be a set of distinct data points or nodes arbitrarily distributed on a domain $ \Omega \subseteq \mathbb{R}^{d}$. Associated with this set is another collection  $ F = \{ f_i = f(\boldsymbol{x}_i) ,i=1, \ldots, N \}$ representing data values obtained by sampling a potentially unknown function  $f: \Omega \rightarrow \mathbb{R}$ at the nodes $ \boldsymbol{x}_i$. The problem at hand is the \emph{scattered data interpolation problem}, which entails discovering an interpolating function $P_f: \Omega \rightarrow \mathbb{R}$ that exactly reproduces the measured values at their respective locations, i.e. 
\begin{equation*}
P_f\left( \boldsymbol{x}_i\right)=f_i, \qquad i=1, \ldots, N.
\end{equation*}

We now suppose to have a univariate function $ \varphi: [0, \infty) \to \mathbb{R}$, known as RBF, which depends on a shape parameter $\varepsilon > 0$ providing, for $\boldsymbol{x},\boldsymbol{z}
\in \Omega$, the real symmetric strictly positive definite kernel 
\begin{equation*} 
\kappa_\varepsilon(\boldsymbol{x},\boldsymbol{z}) = \varphi(\varepsilon ||\boldsymbol{x}-\boldsymbol{z}||_2): = \varphi(\varepsilon r).
\end{equation*}

The kernel-based interpolant $P_f$ can be written as
\begin{equation*}
P_f\left( \boldsymbol{x}\right)= \sum_{k=1}^{N} c_k \kappa_\varepsilon \left( \boldsymbol{x} , \boldsymbol{x}_k  \right), \quad \boldsymbol{x} \in \Omega,
\end{equation*}
whose coefficients $c_k$ are the solution of the linear system
\begin{equation} \label{linsys}
\mK \boldsymbol{c} = \boldsymbol{f},
\end{equation}
where $  \boldsymbol{c}= \left(c_1, \ldots,
c_N\right)^{\intercal}$, $  \boldsymbol{f} =\left(f_1, \ldots , f_N\right)^{\intercal}$, and $\mK_{ik}= \kappa_\varepsilon \left( \boldsymbol{x}_i , \boldsymbol{x}_k  \right)$, $i,k=1, \ldots, N$. Since $\kappa_\varepsilon$ is a symmetric and strictly positive definite kernel, the system \eqref{linsys} has exactly one solution \cite{Fasshauer15}.
Furthermore, the kernel $\kappa_\varepsilon$ gives rise to what is known as the \emph{native space}. This native space, denoted as ${\cal N}{\kappa_\varepsilon}(\Omega)$, is a Hilbert space equipped with the inner product $(\cdot,\cdot)_{{\cal N}{\kappa_\varepsilon}(\Omega)}$. In this space, the kernel $\kappa_\varepsilon$ is reproducing, i.e. for any $f \in {\cal N}_{\kappa_\varepsilon}(\Omega)$ the following identity holds: $f(\boldsymbol{x}) = (f,\kappa_\varepsilon(\cdot,\boldsymbol{x}))_{{\cal N}_{\kappa_\varepsilon}(\Omega)}$, with $\boldsymbol{x}\in \Omega$. By introducing a pre-Hilbert space $H_{\kappa_\varepsilon}(\Omega)= \mbox{span}\{\kappa_\varepsilon(\cdot,\boldsymbol{x}),$ $\boldsymbol{x} \in \Omega\}$, with reproducing kernel $\kappa_\varepsilon$ and equipped with the bilinear form $(\cdot,\cdot)_{H_{\kappa_\varepsilon}(\Omega)}$, the native space ${\cal N}_{\kappa_\varepsilon}(\Omega)$ of $\kappa_\varepsilon$ coincides with its completion with respect to the norm $||\cdot||_{H_{\kappa_\varepsilon}(\Omega)}=\sqrt{(\cdot,\cdot)_{H_{\kappa_\varepsilon}(\Omega)}}$, and for all $f \in {H_{\kappa_\varepsilon}(\Omega)}$ we have $||f||_{{\cal N}_{\kappa_\varepsilon}(\Omega)}=||f||_{H_{\kappa_\varepsilon}(\Omega)}$. 

\subsection{The PUM Scheme}

It is a widely recognized fact that performing the inversion of the kernel interpolation matrix in  \eqref{linsys} can become computationally demanding as the amount of data significantly grows. To address this challenge effectively, a practical approach is to divide the open and bounded domain $\Omega$ into $m$ overlapping subdomains denoted as $\Omega_j$, with the property that $\Omega \subseteq \bigcup_{j=1}^m\Omega_j$. Consequently, this allows for the problem of interpolation to be split independently within each of these subdomains.

The PU covering consists of overlapping balls of radius $\delta$ whose centres are the grid poimts $P=\{\tilde{\boldsymbol{x}}_k$, $k= 1,\ldots,m\}$. In \cite{Fasshauer} it is shown that when the nodes are nearly uniformed distributed, $m$ is a suitable number of PU subdomains on $\Omega$ if ${N}/{m} \approx 2^d$. Then, the covering property is satisfied by taking the radius $\delta$ such that
\begin{equation*} 
      	\delta \geq \frac{\displaystyle 1 }{\displaystyle m^{1/d}}.
      	\label{PU_radius}
\end{equation*}

The PUM solves a local interpolation problem on each subdomain and constructs the global approximant by gluing together the local contributions using weights. To achieve that, we need those weights to be a family of compactly supported, non-negative, continuous functions $w_j$, with $\text{supp}\left(w_j\right) \subseteq \Omega_j$, such that
\begin{equation*}
\sum_{j=1}^{m} w_j\left(\boldsymbol{x}\right) = 1, \quad \boldsymbol{x} \in \Omega.
\end{equation*}
Once we choose the partition of unity $ \{ w_j \}_{j=1}^{m}$, the global interpolant is formed by the weighted sum of $m$ local approximants $P_f^j$, i.e.
\begin{equation*}
P_f\left( \boldsymbol{x}\right)= \sum_{j=1}^{m} P_f^j\left( \boldsymbol{x} \right) w_j \left( \boldsymbol{x}\right) = \sum_{j=1}^{m} \left( \sum_{k=1}^{N_j} c_k^j \kappa_{\varepsilon,\delta}(\boldsymbol{x} ,  \boldsymbol{x}^j_k) \right)  w_j \left( \boldsymbol{x}\right), \quad \boldsymbol{x} \in \Omega,
\label{intg}
\end{equation*}
where $\boldsymbol{x}^j_k \in X_j = X \cap \Omega_j$ and $N_j=|X_j|$, with $k=1, \ldots, N_j$.
We will use the following and well-known Shepard weights in the implementation of our algorithms:
\begin{equation}\label{shep}
    w_j(x)=\frac{\varphi_j(x)}{\sum_{k=1}^m \varphi_k(x)} \quad j=1,\dots,m,
\end{equation}
where $\varphi_k(x)$ is a compactly supported function with support on $\Omega_j$, see \cite{acd,shep68}.

It is worth noting that the accuracy of the fit strongly depends on the choices of the shape parameter and the radius, see e.g. \cite{cav21,cav22,Fornberg-Wright04,Larsson-Fornberg05,lin22}.



An advantage of this scheme is the use of a continuous search in the parameters space $\mathcal{X}= I \times J =(0, \varepsilon_{max}] \times [\delta_{min}, 2 \delta_{min}]$ to obtain a better approximation of $(\varepsilon,\delta)_j^*$, in each subdomain $\Omega_j$.

\section{Bayesian Optimization}
\label{bo}

When seeking to locate a global maximiser for an unknown or challenging-to-assess function $g$ within a bounded set $X$, Bayesian optimization offers an effective approach \cite{Mockus_1978}. Highly regarded in the realm of machine learning, BO is an iterative methodology that optimally utilizes available resources. It entails constructing a probabilistic model of $g$, often referred to as a {\sl surrogate model}, and employing it to guide the selection of sampling points within the set $X$ using an acquisition function. These selected points are located in the area in which the target function will be assessed. After each iteration, the distribution is updated to reflect the acquired information and is subsequently utilized in the next iteration. While some computational effort is required to determine the next point for evaluation, this cost is justifiable when the evaluations of $g$ are computationally expensive. This is because such computations are driven by the goal of reaching the maximum value in a limited number of iterations, which is particularly important in scenarios like optimizing the error function of resource-intensive machine learning algorithms in multi-layer neural networks.

Hereinafter, we briefly review the BO technique \cite{Brochu}. A Gaussian Process (GP) is a collection of random variables such that any subsets of these have a joint Gaussian distribution. Then GPs are completely specified by a mean function $\mathtt{m} : \mathcal{X} \rightarrow \mathbb{R}$ and a positive definite covariance function $k : \mathcal{X} \times \mathcal{X} \rightarrow \mathbb{R}$ (see \cite{Rasmussen}). Further, they are the most common choice for the surrogate model for BO due to the low evaluation cost and the ability to incorporate prior beliefs about the objective function. When modeling the target function with a GP as $g(\textbf{x}) \sim \mathcal{GP} \big( \mathtt{m}(\textbf{x}), k(\textbf{x}, \textbf{x}') \big)$, we impose that
\begin{align*}
    \mathbb{E}\big [ g(\textbf{x}) \big] = \mathtt{m}(\textbf{x}), \qquad
    \mathbb{E}\big [ \big( g(\textbf{x}) - \mathtt{m}(\textbf{x}) \big)  \big( g(\textbf{x}') - \mathtt{m}(\textbf{x}') \big) \big] = k(\textbf{x}, \textbf{x}').
\end{align*}
In the matter of making a prediction given by some observations, the assumption of joint Gaussianity allows retrieving the prediction using the standard formula for mean and variance of a conditional normal distribution. Hence, suppose to have $s$ observation $\boldsymbol{g} = (g(\textbf{x}_1), \dots, g(\textbf{x}_s))^{\intercal}$ on $\textbf{X} = (\textbf{x}_1, \dots, \textbf{x}_s)^{\intercal}$ and a new point $\bar{\textbf{x}}$ on which we are interested in having a prediction of $\bar{g} = g(\bar{\textbf{x}})$. The previous observations $\boldsymbol{g}$ and the predicted value $g(\bar{\textbf{x}})$ are jointly normally distributed:
$$
Pr\begin{pmatrix} \begin{bmatrix} \boldsymbol{g} \\ g(\bar{\textbf{x}}) \end{bmatrix} \end{pmatrix} = \mathcal{N} \begin{bmatrix}
\begin{bmatrix} \mu(\textbf{X}) \\ \mu(\bar{\textbf{x}}) \end{bmatrix}, 
\begin{bmatrix} 
K(\textbf{X},\textbf{X}) \ \ K(\textbf{X}, \bar{\textbf{x}}) \\
K(\textbf{X}, \bar{\textbf{x}})^{\intercal} \ \ k(\bar{\textbf{x}}, \bar{\textbf{x}}) \\
\end{bmatrix}\end{bmatrix},
$$
where $K(\textbf{X},\textbf{X})$ is the $s \times s$ matrix with $(i, j)$-element $k(\textbf{x}_i, \textbf{x}_j)$, and
$K(\textbf{X}, \bar{\textbf{x}})$ is a $s \times 1$ vector whose $i^{\text{th}}$ element is given by $k(\textbf{x}_i, \bar{\textbf{x}})$, see \cite{Rasmussen}. Since $Pr(f(\bar{\textbf{x}}) |  \boldsymbol{g})$ must also be normal, it is also possible to estimate the distribution, the mean and the covariance, for any point in the domain. When data points and data values retrieved by the evaluation of the target function are fed to the model, they induce a posterior distribution over functions which is used for the next iteration as a prior. It is worth noting that in the case of modeling a function with a GP, when we observe a value, we are essentially observing the random variable associated with that specific point.



An \emph{acquisition function} $a: \mathcal{X} \rightarrow \mathbb{R}$ serves as a tool for determining the subsequent point at which the objective function will be assessed. The goal is to choose a point that maximizes this acquisition function, and the result of evaluating the objective function at this chosen point is utilized to update the surrogate model. The design of an acquisition function is specifically crafted so that a high acquisition score corresponds to the likelihood of encountering high values of the objective function. When the decision is made regarding which acquisition function to use, a trade-off arises between exploration and exploitation. Exploration involves the selection of points characterized by high levels of uncertainty, typically those located at a considerable distance from previously examined points. Conversely, exploitation involves the selection of points in close proximity to those already assessed by the objective function. The most common acquisition functions are:

\begin{itemize}
    \item \textbf{Probability of Improvement}, which maximises the probability of improvement over the best current value;
    \item \textbf{Expected Improvement}, which maximises the expected improvement over the current best;
    \item \textbf{GP Upper Confidence Bound}, which minimises the cumulative regret\footnote[1]{Regret is a performance metric commonly used in Reinforcement Learning. In a maximization setting of a function $g$ it represents the loss in rewards due to not knowing $g$'s maximum points beforehand. If $x^{\ast} = \argmax g(x)$, the regret for a point $x$ is $g(x^{\ast}) - g(x)$ over the course of the optimization.}.
\end{itemize}
In what follows we will use the \lq\lq Expected Improvement\rq\rq\ \cite{EI} as the acquisition function. Not only does it consider the probability of improvement of the candidate point with respect to the previous maximum, but also the magnitude of this improvement.

Suppose that after a number of iterations the current maximum of the objective function is $g(\hat{\textbf{x}})$. Given a new point $\textbf{x}$, the Expected Improvement acquisition function computes the expectation of improvement $g(\textbf{x}) - g(\hat{\textbf{x}})$  over the part of the normal distribution that is above the current maximum:
\begin{equation}
\label{EI_eq}
    EI(\textbf{x}) = \int_{g(\hat{\textbf{x}})}^\infty \big( g^*(\textbf{x}) - g(\hat{\textbf{x}}) \big) 
    \frac{1}{\sqrt{2 \pi} \sigma(\textbf{x})} e^{-\frac{1}{2}  [(g^*(\textbf{x}) - \mu(\textbf{x}))/\sigma(\textbf{x})]^2} dg^*(\textbf{x}),
\end{equation}
where $g^*(\textbf{x})$, $\mu(\textbf{x})$ and $\sigma(\textbf{x})$ represent the predicted value by the surrogate model, the expected value and the variance of $\textbf{x}$, respectively.
Solving integral \eqref{EI_eq} leads to the following closed form for the evaluation of the Expected Improvement:
\begin{equation*} \label{EI_closed}
    EI(\textbf{x}) = 
    \begin{cases}  
            (\mu(\textbf{x}) - g(\hat{\textbf{x}})) \Phi(Z) + \sigma(\textbf{x)} \phi(Z), & \text{ if } \sigma(\textbf{x}) >0, \\
            0, & \text{ if } \sigma(\textbf{x}) =0, 
    \end{cases}
\end{equation*}
where $Z =\frac{\mu(\textbf{x})-g(\hat{\textbf{x}})}{\sigma(\textbf{x})}$, while $\phi$ and $\Phi$ are the Probability Density Function and Cumulative Distribution Function of the standard normal distribution $\mathcal{N}(0, 1)$, respectively. 
An extension of (\ref{EI_closed}) that also trades off exploration and expectation by means of a non-negative parameter $\xi$ was proposed in \cite{Lizotte}: 
\begin{equation*} 
\label{EI_closed_trade-off}
   EI(\textbf{x}) = 
   \begin{cases}  
           (\mu(\textbf{x}) -g(\hat{\textbf{x}}) -\xi) \Phi(Z) + \sigma(\textbf{x}) \phi(Z), & \text{ if } \sigma(\textbf{x}) >0, \\
           0, & \text{ if } \sigma(\textbf{x}) =0, 
   \end{cases}
\end{equation*}
where $Z =\frac{\mu(\textbf{x})-g(\hat{\textbf{x}}) -\xi}{\sigma(\textbf{x})}$.

\section{Algorithms and Their Computational Cost}
\label{algorithmsandcc}


In this section, we first describe in Subsection \ref{algorithms} the algorithms for interpolation processes with Bayesian optimization. In Subsection \ref{cc} their computational cost is analysed.

\subsection{Algorithms}\label{algorithms}
Let $X$ be a set of points for which we know the associated set of data values $F$, and let $\bar{X} = \{  \boldsymbol{\bar{x}}_i, i = 1,  \ldots , \bar{n} \}$ be a set of points on which we want to evaluate some approximate solutions.
The whole process is handled by Algorithm \ref{bo_pum}, which invokes the Bayesian optimization (Algorithm \ref{bo_alg}) for the parameters search and the partition of unity (Algorithm \ref{pum}) for the evaluation of the approximant. In detail, Algorithm \ref{bo_pum} builds the approximant on $\bar{X}$ with the best shape parameters $\bm{\varepsilon}$ and radii $\bm{\delta}$ found by means of the Bayesian optimization. With the aim of doing that, the algorithm starts by retrieving $N$, the number of points in $X$, and the dimension $d$ of the space. Using these values, it evaluates the number $m$ of partition of unity centers and generates them as an equally spaced grid in $\Omega$. For the sake of clarity, hereinafter, without loss of generality we consider the special case of $\Omega=[0,1]^d$.  We remark that a suitable number of subdomains is $\lfloor \frac{N}{2^{d}}\rfloor$, see \cite{Fasshauer}.
After that, it evaluates the distance tree of $X$, the function \emph{KDTree} of the package \emph{scipy.spatial} and its method are used to build and perform the points search on it. Next, in each subdomain the value for the radius that ensures the minimum density is found by Algorithm \ref{find_min_radius}. This value $\bm{\delta}_{start}$ will be used when applying a Bayesian optimization (Algorithm \ref{bo_alg}) to enhance the shape parameter and the radius in each subdomain. The last step is to train the RBF-PUM approximant (Algorithm \ref{rbf}) with the found parameters and return the approximated value of the function on the set of points $\Bar{X}$.


\begin{algorithm}
\caption{BO-PUM}
\hspace*{\algorithmicindent} \textbf{Input:} 
    \vspace{-3mm}
\begin{itemize}
\setlength\itemsep{-0.2em}
  \item[] $X$: data points,
  $F$: data values,
   $\Bar{X}$: evaluation points, $I$: $\varepsilon$ search interval,   $a$: acquisition function, $\xi$: exploration-exploitation parameter, $nstart$: number of starting points,  $niter$: number of Bayesian iterations, $min_{pts}$: minimum number of points in a subdomain, $\tau$: tolerance of the error during the parameters search, $w$: weight function. 
\end{itemize}

\begin{algorithmic}
\State $N \rightarrow$ number of points in $X$ 
\State $d \rightarrow$ space dimension of $X$
\State $m \rightarrow \lfloor \frac{N}{2^{d}}\rfloor$
\State $centers \rightarrow$ grid of $m$ points in $\Omega= [0,1]^d$
\State $\bm{\varepsilon} \rightarrow$ $0$-vector of length $m$
\State $\bm{\delta} \rightarrow$ $0$-vector of length $m$
\State $DT_X \rightarrow$ distance tree of $X$
\State $\bm{\delta}_{start} \rightarrow \textbf{FIND-MIN-RADIUS}(X, centers,m,d,min_{pts},DT_{{X}})$ 
\For{$i = 1 : |centers|$}
    \State $J \rightarrow [\bm{\delta}_{start}[i], 2 \bm{\delta}_{start}[i]]$
    \State $[\varepsilon, \delta] \rightarrow \textbf{BO}(X, F,I,J,a,\xi,nstart, niter, DT_{X}, centers[i], \tau)$
    
    \State $\bm{\varepsilon}[i] \rightarrow \varepsilon$
    \State $\bm{\sigma}[i] \rightarrow \sigma$
\EndFor
\State $\textbf{P}_{\boldsymbol{f}} \rightarrow \textbf{PUM}(X, F,\Bar{X}, centers, DT_X, w, \bm{\varepsilon}, \bm{\delta})$

\end{algorithmic} 
\hspace*{\algorithmicindent} \textbf{Output:} 
    \vspace{-3mm}
\begin{itemize}
\setlength\itemsep{-0.2em}
    \item[] $\textbf{P}_{\boldsymbol{f}}$: evaluation of the interpolated solution on $\Bar{X}$
\end{itemize}
\label{bo_pum}   
\end{algorithm}

\begin{algorithm}
\caption{FIND-MIN-RADIUS}
\hspace*{\algorithmicindent} \textbf{Input:} 
    \vspace{-3mm}
\begin{itemize}
\setlength\itemsep{-0.2em}
    \item[] $X$: data points,  $centers$: PU centers, $m$: number of centers,  $d$: space dimension, $min_{pts}$: minimum number of points in a subdomain, $DT_{X}$: distance tree of $X$.
\end{itemize}
\begin{algorithmic}
\State $\bm{\delta}_{start} \rightarrow$ $\frac{\sqrt{d}}{2m^\frac{1}{d}} \times 1$-vector of length $|centers|$
    \For{$i = 1 : m$}
    \State $X_{sub} \rightarrow$ retrieve the subset of ${X}$ within distance $\bm{\delta}_{start}[i]$ from $centers[i]$ using $DT_{{X}}$
    \While{$|X_{sub}| < min_{pts}$}
        \State $\bm{\delta}_{start}[i] \rightarrow \bm{\delta}_{start}[i] + \frac{1}{8} \frac{\sqrt{d}}{2m}$
        \State $X_{sub} \rightarrow$ retrieve the subset of ${X}$ within distance $\bm{\delta}_{start}[i]$ from $centers[i]$ using $DT_{{X}}$
    \EndWhile
\EndFor
\end{algorithmic}
\hspace*{\algorithmicindent} \textbf{Output:} 
    \vspace{-3mm}
\begin{itemize}
\setlength\itemsep{-0.2em}
    \item[] $\bm{\delta}_{start}$: vector of subdomain radii that ensure the minimum densities. 
\end{itemize}
\label{find_min_radius}
\end{algorithm}

The core of the process is accomplished by Algorithm \ref{bo_alg}, which is a remodelling of the  \emph{BayesianOptimization} Python's library \cite{BOPy}. Concretely, it traces the optimisation process provided by the \emph{optimisation} method of the \emph{BayesianOptimisation} class, which in sequence uses the methods \emph{fit} and \emph{predict} of the function \emph{GaussianProcessRegressor} of the \emph{sklearn.gaussian\_process} package \cite{sklearn}. Further details about the implementation of \emph{GaussianProcessRegressor} are available at \cite[Algorithm 2.1]{Rasmussen}.

To measure the goodness of the approximant, we introduce the Maximum Absolute Error (MAE), the Relative Maximum Absolute Error (RMAE) and the Relative Root Mean Squared Error (RRMSE) defined as follows: 
\begin{equation*}
    \mbox{MAE}(X, F, \textbf{P}_{\boldsymbol{f}} ) = \mbox{MAE}_{X, F}(\textbf{P}_{\boldsymbol{f}}) = \max_{\boldsymbol{x}_i \in X, f_i \in F} | P_f(\boldsymbol{x}_i)- f_i|,
\end{equation*} 
\begin{equation*}
    \mbox{RMAE}(X, F, \textbf{P}_{\boldsymbol{f}} ) = \mbox{RMAE}_{X, F}(\textbf{P}_{\boldsymbol{f}}) = \max_{\boldsymbol{x}_i \in X, f_i \in F} \frac{| P_f(\boldsymbol{x}_i)- f_i|}{f_i},
\end{equation*} 
\begin{equation*}
    \mbox{RRMSE}(X, F, \textbf{P}_{\boldsymbol{f}} ) = \mbox{RRMSE}_{X, F}(\textbf{P}_{\boldsymbol{f}}) = \sqrt{\frac{1}{N}\sum_{i=1}^{N}  \bigg(\frac{P_f(\boldsymbol{x}_i)- f_i}{f_i}\bigg)^2},
\end{equation*} 
where $X$ and $F$ are the given sets of data points and data values and  $\textbf{P}_{\boldsymbol{f}} = (P_f(\boldsymbol{x}_1), \dots P_f(\boldsymbol{x}_N))$, with $N = |X|$.
The first step is to initialize the function $g$ to optimize, defined as the negative $\mbox{MAE}_{X_{val}, F_{val}}(\cdot)$ between the known values $F_{val}$ and the approximation evaluated by Algorithm \ref{rbf} for the points $X_{val}$ for a specific value of $\theta$, and the search space $\mathcal{X}$. We remark that $X_{val}, F_{val}$ are subsets of $X$ used for the BO search. Then, until the number of iterations is reached or the error drops below a certain tolerance, for the first $nstart$ iteration the algorithm randomly samples $\hat{\theta}=(\hat{\varepsilon}, \hat{\delta})$  in the search domain $\mathcal{X}$, otherwise the chosen point $\hat{\theta}$ is the one that maximises the acquisition function evaluated on a random set (the acquisition function exploits the Gaussian process fitted in the previous iteration). At this point the algorithm retrieves the subsets $X_j \subset X$ and $ \tilde{X}_j \subset \tilde{X}$ that are contained in the related subdomain $\Omega_j$, splits them into training and validation subsets and fits an approximant applying Algorithm \ref{rbf}.
At each iteration the values of $\hat{\theta}$ and the error obtained fitting the interpolant with the parameter $\hat{\theta}$ are stored in the vectors $\bm{\theta}$ and $\bm{g}$ and a Gaussian process on $(\bm{\theta}, \bm{g})$ is fitted.
The last step consists of determining $\theta^*$, the parameter that maximises the vector $\bm{g}$.

\begin{algorithm}
\caption{BO}
\hspace*{\algorithmicindent} \textbf{Input:} 
    \vspace{-3mm}
\begin{itemize}
\setlength\itemsep{-0.2em}
    \item[] $X$: data points, $F$: data values,  $I$: parameter search interval for $\varepsilon$, $J$: parameter search interval for $\delta$, $a$: acquisition function, $\xi$: exploration-exploitation parameter, $nstart$: number of starting points, $niter$: number of Bayesian iterations, $DT_{X}$: distance tree of $X$,  $center$: subdomain center, $\tau$: tolerance.
\end{itemize}
\begin{algorithmic}

\State Set $g \rightarrow -\mbox{MAE}_{X_{val}, F_{val}}(\cdot)$
\State $\mathcal{X} \rightarrow  I \times J$
\State $\bm{g} \rightarrow (\cdot)$ (empty vector)
\State  $\bm{\theta} \rightarrow (\cdot)$ (empty vector)
\While{$i \leq nstart + niter \text{\qquad or \qquad} |\max(\bm{g})| > \tau$}
    \If{$i \leq nstart$}
        \State $\hat{\theta} \rightarrow$ random sample in $\mathcal{X}$ (note that $\hat{\theta} = (\hat{\varepsilon}, \hat{\delta})$)
    \Else
        \State Evaluate $a$ on a set of random points in $\mathcal{X}$
        \State Select the point $\hat{\theta}$ that maximizes $a$
    \EndIf    
    \State $X_{sub} \rightarrow$ retrieve the subset of $X$ within distance $\hat{\delta}$ from $center$ using $DT_{X}$
    \State Split $X_{sub}, F_{sub}$ in $X_{train},  X_{val},  F_{train}, F_{val}$
    \State $\textbf{P}_{\boldsymbol{f}_{\hat{\theta}}}  \rightarrow \textbf{RBF}(X_{train}, F_{train}, X_{val}, \hat{\varepsilon})$  \qquad (call to \textbf{Algorithm} \ref{rbf})

    \State $\bm{\theta} \rightarrow \bm{\theta} \cup \hat{\theta} $
    \State $\bm{g} \rightarrow \bm{g} \cup g({\textbf{P}_{\boldsymbol{f}_{\hat{\theta}}}})) $
    \State Fit the Gaussian process on $(\bm{\theta}, \bm{g})$
    \State $ i \rightarrow i + 1$
\EndWhile
\State $\theta^* \rightarrow \argmax \bm{g}$
\end{algorithmic}
\hspace*{\algorithmicindent} \textbf{Output:} 
    \vspace{-3mm}
\begin{itemize}
\setlength\itemsep{-0.2em}
    \item[] $\theta^*$: best parameters.
\end{itemize}
\label{bo_alg}
\end{algorithm}

The fundamental component of the scheme is represented by Algorithm \ref{rbf}, which solves the interpolation system, and finds the approximated values of the function on $\bar{X}$.

\begin{algorithm}
    \caption{RBF}
    \hspace*{\algorithmicindent} \textbf{Input:}
    \vspace{-3mm}
    \begin{itemize}
        \setlength\itemsep{-0.2em}
        \item[] $X$: data points, $F$: data values, $\Bar{X}$: evaluation points,  $\varepsilon$: shape parameter.   
    \end{itemize}
    
    \begin{algorithmic}
        \State  $\textbf{P}_{\boldsymbol{f}} \rightarrow$ Solve the linear system of the form \eqref{linsys}
    \end{algorithmic}
    \hspace*{\algorithmicindent} \textbf{Output:}
        \vspace{-3mm}
\begin{itemize}
\setlength\itemsep{-0.2em}
    \item[] $\textbf{P}_{\boldsymbol{f}}$: evaluation of the interpolated solution onapproximate evaluation on $\Bar{X}$.
    \end{itemize}
\label{rbf}
\end{algorithm}


Algorithm \ref{pum}, given the values of $\varepsilon$ and $\delta$ for each subdomain, deals with determining local solutions and summing them up to construct a global one. In detail, it determines the Shepard weights in \eqref{shep} for the subdomains and the approximate evaluation array on $\bar{X}$ is initialised with all zeros.

It is worth remarking that the algorithms presented in this section can be extended to manage the approximation settings by selecting a subset of $X$ as the set of centres for the RBF and solving the linear system (\ref{linsys}) in the least squares setting. We refer the reader to \cite{Fasshauer15} for further details.  

\begin{algorithm}
\caption{PUM}
\hspace*{\algorithmicindent} \textbf{Input:} 
    \vspace{-3mm}
\begin{itemize}
\setlength\itemsep{-0.2em}
  \item[] $X$: data points, 
  $F$: data values,
   $\Bar{X}$: evaluation points, $centers$: subdomain centers, $DT_X$: distance tree of $X$,  $w$: weight function, $\bm{\varepsilon}$: vector of shape parameters,  $\bm{\delta}$: vector of subdomain radius.
\end{itemize}

\begin{algorithmic}
\State $\textbf{sw} \rightarrow$  retrieve Shepard weights from $w$
\State $\textbf{P}_{\boldsymbol{f}} \rightarrow$ $0$-vector of length $|\bar{X}|$
\State $DT_{\bar{X}} \rightarrow$ distance tree of $\bar{X}$
\For{$j = 1:|centers|$}
    \State $n_j \rightarrow$ indices of points of $X$ at distance $\bm{\delta}[j]$ from $centers[j]$  using $DT_{X}$

\If{$|n_j| \neq 0$} 
    \State $s_j \rightarrow$ points of $\bar{X}$ at distance $\bm{\delta}[j]$ from $centers[j]$ using $DT_{\bar{X}}$
        \If {$(|s_j| \neq 0)$}
            \State $\textbf{P}_{\boldsymbol{f}_{\textbf{s}_\textbf{j}}}' \rightarrow \textbf{RBF}(X_{j},  F_{j}, \bar{X}_{s_j}, \varepsilon_j )$ 
            \State $\textbf{P}_{\boldsymbol{f}_{\textbf{s}_\textbf{j}}} \rightarrow \textbf{P}_{\boldsymbol{f}_{\textbf{s}_\textbf{j}}} + \textbf{P}_{\boldsymbol{f}_{\textbf{s}_\textbf{j}}}' * \textbf{sw}_{\textbf{s}_\textbf{j}}$
        \EndIf
\EndIf
\EndFor

\end{algorithmic} 
\hspace*{\algorithmicindent} \textbf{Output:} 
    \vspace{-3mm}
\begin{itemize}
\setlength\itemsep{-0.2em}
    \item[] $\textbf{P}_{\boldsymbol{f}}$: evaluation of the interpolated solution on $\Bar{X}$
\end{itemize}
\label{pum}
\end{algorithm}

\subsection{Computational Analysis of Algorithms}\label{cc}
In this section we discuss the complexity of the Algorithm \ref{bo_pum} by first analyzing its constituent components. We remark here that logarithms are taken base 2. Without explicitly declaring the shape parameter and the subdomain radius, we will show that the total expense required to build a global interpolant is $\mathcal{O}(N^\frac{2d + 1}{d} + N(nstart +niter)(N+N_j^ 3 + (nstart+niter)^3))$, where $N = |X|$, $d$ is the space dimension, $N_j$ is the maximum among the number of points in the subdomains (note that $N_j << N$), and  $nstart$ and $niter$ are the initialization and the Bayesian steps of the optimization.  We will need to consider the cost of construction and search in a kdtree, $\mathcal{O}(N \log(N))$ and $\mathcal{O}(N)$, respectively. We use a kdtree implementation by \emph{scipy} \cite{scipy} that provides, for a balanced dataset, a balanced tree by applying a median-based splitting strategy in $\mathcal{O}(N \log(N))$. We want to highlight that hardly ever and very uncommon in practice, the computational expenses for the construction could be $\mathcal{O}(N^2)$ in the worst case. A different discussion can be addressed in the case of search, where we face three different scenarios: the best-case scenario, where for balanced tree, the search cost $\log(N)$; the average-case, where the kdtree is reasonably balanced and the search cost is $\mathcal{O}(N + K)$ where $K$ is the number of points found in the search distance; the worst-case, where the tree is unbalanced and the computational cost is $\mathcal{O}(N)$.\\

\textbf{FIND-MIN-RADIUS:}
Suppose that the initial radius for a subdomain is equal to $0$. In this case, the maximum number of augmenting steps inside the $while$ loop is bounded by the length of the diagonal of the $d$-dimensional hypercube over the weight of the extent. Hence it is bounded with $16 m^\frac{1}{d} \simeq 8 N^\frac{1}{d}$.
Taking into account that the search has a cost of $N$, we can estimate the cost of the whole radius search in Algorithm \ref{find_min_radius} as: 

\begin{align*}
\begin{array}{rcl}
    \mathcal{O}(m + m[N + 16 m^\frac{1}{d} N]) & \simeq  & \mathcal{O}( m[N + 16 m^\frac{1}{d} N]) \medskip \\
                                                & \simeq & \mathcal{O}( m^\frac{d + 1}{d} N) \medskip \\
                                                & \simeq & \mathcal{O}( N^\frac{2d + 1}{d}). 
\end{array}
\end{align*}

\textbf{RBF:}
Algorithm \ref{rbf} simply involves the solution of a linear system. With an input of $N$ nodes the computational expense is $\mathcal{O}(N^3)$.\\

\textbf{PUM:}
Overlooking the cost of the computation of the Shepard weights and some initialization that has cost linearly dependent from $N$, we have a $for$ loop of length $m$ where for each iteration we have three point search of cost $N$, a call of Algorithm \ref{rbf} with a variable input dimension and an update of a $N$-length vector. Let $\tilde{N}= \max_j N_j = \max_j |X_j|$ be, i.e. it is  the maximum among the number of points in the subdomain. We have that the complexity of Algorithm \ref{pum} is:
\begin{align*}
\begin{array}{rcl}
    \mathcal{O} ( m[N + \tilde{N}^3]) & \simeq &  \mathcal{O} ( N[N + \tilde{N}^3]) \medskip\\
                                & \simeq &  \mathcal{O} ( N^2 + N\tilde{N}^3]).
\end{array}
\end{align*}

\textbf{BO:}
Algorithm \ref{bo_alg} is made by a $while$ loop of length less than $nstart + niter$ iteration. For each iteration the cost can be summarised as follows: computation of order $N$ for selecting the next parameter to evaluate,  $2$ points search, a splitting of cost $N$, an invocation of Algorithm \ref{rbf} of cost $\tilde{N}^3$ and the fitting of the Gaussian process that in the worst case cost $(nstart + niter)^3$.
The total expense for the Bayesian optimization is:
\begin{align*}
    \mathcal{O}((nstart + niter)(N + \tilde{N} + (nstart+niter)^3)).
\end{align*}

\textbf{BO-PUM:}
In conclusion, summing up the previous results, and taking into account that the construction of the kdtrees requires computation of order $N \log(N)$,  we can retrieve the computational expense for Algorithm \ref{bo_pum} as follows:
\begin{align*}
   & \mathcal{O}( N + N \log(N) + N^\frac{2d+1}{d} + m (start + niter)(N + \tilde{N}^3 + (nstart+niter)^3) + N^2 + N\tilde{N}^3)   \\
   & \qquad \qquad  \simeq \mathcal{O}(N^\frac{2d+1}{d} + N (start + niter)(N + \tilde{N}^3 + (nstart+niter)^3)  \\
   & \qquad \qquad  \simeq \mathcal{O}(N^\frac{2d+1}{d} + N^2 (start + niter) +N\tilde{N}^3 (start + niter) + N(nstart+niter)^4).
\end{align*}




\section{Numerical Experiments and Applications} \label{exp}
In this section, we will illustrate how algorithms presented in Section \ref{algorithmsandcc} work efficiently both on test and real-world datasets. 

Before going into details, we remark that all the code was developed in Python 3.9 and the library used to perform the optimization is {\em BayesianOptimization} \cite{BOPy} in which the default kernel used for the Gaussian process is the Matérn $5/2$. Moreover, we set the parameters $\xi=0.15$ and $min_{pts}$= 15. To apply BO in the search for optimal parameters $(\varepsilon, \delta)$, we assume that the objective function to be maximized is the Maximum Absolute Error (MAE) of the RBF interpolant, with the sign inverted. This is because BO is a maximization process, as explained in Section \ref{bo}. We vary the number of points in the training set while keeping the test set fixed at 1000 points. During the BO process, for each subdomain, after identifying the points within it, we further divide them into sub-training and sub-validation sets to enable the evaluation of the training error. After determining the best parameter pairs for each subdomain, we train a PUM interpolant on the training set for each optimizer using the identified parameters.
For each subdomain, the search space is $ \mathcal{X}=(0, 20] \times [\delta_{min}, 2 \delta_{min}]$, where $\delta_{min}$ is the radius value that ensures a minimum density in the subdomain. 
BO performs $5$ random steps plus at most $25$ Bayesian steps in the search space. The iterative process stops when the desired tolerance $\tau$ is reached.
\subsection{Numerical Experiments}
We perform the experiments on four different sizes of random data in the domain $\Omega=[0, 1]^2$ using three RBFs of different smoothness, i.e.,
\begin{align*}
&\varphi_1(\varepsilon r) = e^{-\varepsilon^2 r^2} & (\mbox{Gaussian $C^{\infty}$}), \\ 
&\varphi_2(\varepsilon r)= e^{-\varepsilon r}(1 + 3 \varepsilon r + \varepsilon^2 r^2) & (\mbox{Mat$\acute{\text{e}}$rn $C^4$}),\\
& \varphi_3(\varepsilon r) = (35 \varepsilon^2 r^2 + 18 \varepsilon r + 3)  |1 - \varepsilon r|_+^6  & (\mbox{Wendland $C^4$}),
\end{align*}
and the following test functions \cite{laz02,ren99}:
\begin{align*}
&   f_1(x_1,x_2)   =   0.75 \exp\left[{-\frac{(9x_1-2)^2}{4}-\frac{(9x_2-2)^2}{4}}\right] +0.75 \exp\left[{-\frac{(9x_1-2)^2}{49} - \frac{9x_2+1}{10}}\right]  \medskip \\ 
& \qquad \qquad \quad  + 0.5 \exp\left[{-\frac{(9x_1-7)^2} {4}-\frac{(9x_2-3)^2}{4}}\right] -0.2 \exp\left[{-(9x_1-4)^2-(9x_2-7)^2}\right], \bigskip  \\
&     f_2(x_1,x_2) =  2 \cos(10 x_1) \sin(10 x_2) + \sin(10 x_1 x_2). 
\end{align*}

Results are shown in Tables \ref{tab:franke}, \ref{tab:f_2} and \ref{tab:W4}. It is worth noting that as the number of points increases, the execution time of the BO decreases. This is due to the high density of the space when a greater number of points is considered.
 In particular, when this happens, there are denser subdomains, and thus better accuracy and fewer BO iterations are needed to satisfy the tolerance $\tau$.

\begin{table}
    \centering
    
\begin{tabular}{|*{6}{p{17mm}|}}

\cline{3-6}
    \multicolumn{2}{c}{} & \multicolumn{2}{|c}{Gaussian kernel ($\varphi_1$)} &\multicolumn{2}{|c|}{Matérn kernel ($\varphi_2$)}\\
    \hline   
$N$ &     $\tau$ &  time (s)  & MAE  & time (s) & MAE \\
\hline
\multirow{2}{*}{2000} &    1e-04 & 1.02e+01 &  8.16e-05 &  5.56e+01 &  2.15e-04 \\
\cline{2-6}
\multirow{2}{*}{} &   1e-05 &   6.92e+01 &  1.00e-05 &  3.13e+02 &  1.66e-04 \\
\hline
\multirow{2}{*}{4000} &  1e-04 &   4.35e+00 &  2.68e-05 &  1.83e+01 &  6.81e-05 \\
\cline{2-6}
\multirow{2}{*}{} &    1e-05 &   2.43e+01 &  5.50e-06 &  4.49e+02 &  4.36e-05 \\
\hline
\multirow{2}{*}{8000} &  1e-04 & 2.87e+00 &  9.14e-06 &  5.21e+00 &  3.28e-05 \\
\cline{2-6}
\multirow{3}{*}{} &    1e-05 &  1.01e+01 &  5.49e-06 &  3.59e+02 &  3.00e-05 \\
\hline
\multirow{2}{*}{16000} &    1e-04 &  5.54e+00 &  1.25e-06 &  6.16e+00 &  3.59e-05 \\
\cline{2-6}
\multirow{3}{*}{} &    1e-05 & 6.31e+00 &  1.07e-06 &  1.07e+02 &  2.07e-05 \\
\hline
\end{tabular}
    \caption{Computational time and MAE using BO optimizer for Gaussian and Matérn kernels and different number $N$ of random points in $\Omega=[0,1]^2$ using $f_1$ test function. Two tolerances $\tau$ for the training error are used.}
    \label{tab:franke}
\end{table}

\begin{table}
    \centering
    
\begin{tabular}{|*{6}{p{17mm}|}}

\cline{3-6}
    \multicolumn{2}{c}{} & \multicolumn{2}{|c}{Gaussian kernel ($\varphi_1$)} &\multicolumn{2}{|c|}{Matérn kernel ($\varphi_2$)}\\
    \hline   
$N$ &     $\tau$ &  time (s)  & MAE  & time (s) & MAE \\
\hline
\multirow{2}{*}{2000} &   1e-04 &    3.79e+01 &  7.14e-05 &  3.98e+02 &  1.84e-02 \\
\cline{2-6}
\multirow{2}{*}{} &   1e-05 &  2.39e+02 &  3.62e-04 &  3.97e+02 &  1.02e-02 \\
\hline
\multirow{2}{*}{4000} &     1e-04 & 1.35e+01 &  3.16e-05 &  6.73e+02 &  2.29e-03 \\
\cline{2-6}
\multirow{2}{*}{} &     1e-05 & 1.32e+02 &  8.83e-06 &  7.55e+02 &  1.53e-03 \\
\hline
\multirow{2}{*}{8000} &   1e-04 &  6.18e+00 &  9.40e-05 &  6.43e+02 &  8.56e-04 \\
\cline{2-6}
\multirow{3}{*}{} &    1e-05 & 5.47e+01 &  9.63e-06 &  1.46e+03 &  8.84e-04 \\
\hline
\multirow{2}{*}{16000} &    1e-04 & 5.65e+00 &  1.09e-05 &  1.36e+02 &  8.06e-05 \\
\cline{2-6}
\multirow{3}{*}{} &    1e-05 & 1.87e+01 &  5.41e-06 &  2.78e+03 &  1.22e-04 \\
\hline
\end{tabular}
    \caption{Computational time and MAE using BO optimizer for Gaussian and Matérn kernels and different number $N$ of random points in $\Omega=[0,1]^2$ using $f_2$ test function. Two tolerances $\tau$ for the training error are used.}
    \label{tab:f_2}
\end{table}

\begin{table}
    \centering
    
\begin{tabular}{|*{6}{p{17mm}|}}

\cline{3-6}
    \multicolumn{2}{c}{} & \multicolumn{2}{|c}{$f_1$} &\multicolumn{2}{|c|}{$f_2$}\\
    \hline   
$N$ &     $\tau$ &  time (s)  & MAE  & time (s) & MAE \\
\hline
\multirow{2}{*}{2000} &   1e-04 &    2.26e+02 &      1.35e-03 &  4.57e+02 &  1.07e-02 \\
\cline{2-6}
\multirow{2}{*}{} &   1e-05 &  4.06e+02 &      1.57e-02 &  4.52e+02 &  3.11e-02 \\
\hline
\multirow{2}{*}{4000} &     1e-04 & 2.47e+02 &      3.49e-03 &  8.55e+02 &  2.81e-03 \\
\cline{2-6}
\multirow{2}{*}{} &     1e-05 & 7.14e+02 &      7.10e-04 &  8.92e+02 &  3.95e-03 \\
\hline
\multirow{2}{*}{8000} &   1e-04 & 2.70e+02 &      4.15e-04 &  1.25e+03 &  1.27e-02 \\
\cline{2-6}
\multirow{3}{*}{} &    1e-05 & 1.05e+03 &      1.23e-03 &  1.74e+03 &  2.41e-03 \\
\hline
\multirow{2}{*}{16000} &    1e-04 & 3.25e+02 &      1.59e-04 &  1.26e+03 &  8.06e-04 \\
\cline{2-6}
\multirow{3}{*}{} &    1e-05 & 1.22e+03 &      1.15e-04 &  3.37e+03 &  7.25e-04 \\
\hline
\end{tabular}
    \caption{Computational time and MAE using BO optimizer for Wendland kernel $\varphi_3$ and different number $N$ of random points in $\Omega=[0,1]^2$ using $f_1$ and $f_2$ test functions. Two tolerances $\tau$ for the training error are used.}
    \label{tab:W4}
\end{table}

\subsection{Real Data Applications} \label{real_data_appl}
In this subsection we show the behaviour of our framework PUM-BO applied on two different real data examples showing the performance of the algorithm when the measurements are taken with regular intervals, similar to grid points, and on contour lines, similar to random measurements.\\

\textbf{Tonga Trench Dataset:} The Tonga Trench, situated within the vast expanse of the Pacific Ocean, descends to an astonishing depth of 10,882 meters (35,702 feet) at its lowest point, aptly named Horizon Deep. This trench, accompanied by an adjacent volcanic island arc, constitutes an active subduction zone nestled between two tectonic plates within Earth's lithosphere.



\begin{figure}
    \centering
    \includegraphics[width=.48\textwidth]{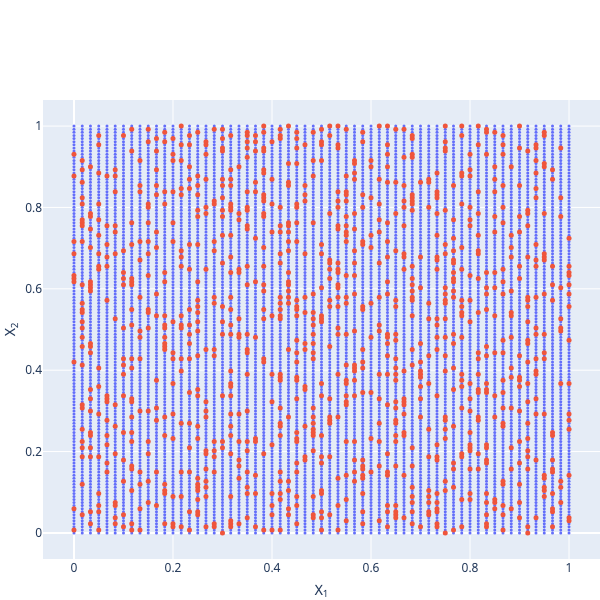} 
\includegraphics[width=.48\textwidth]{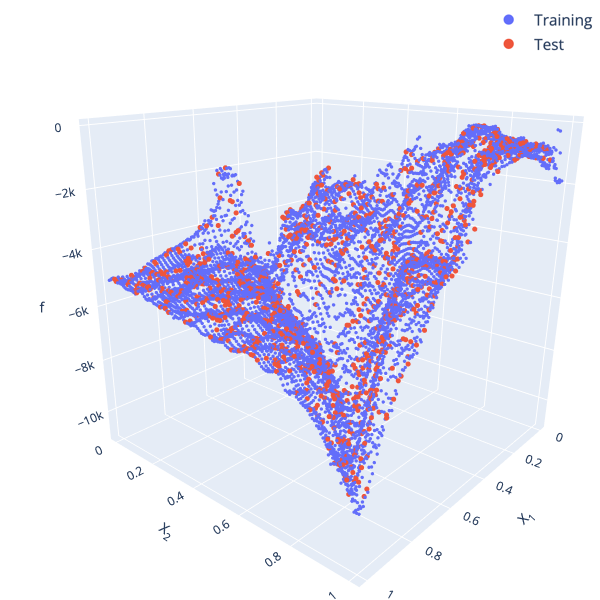}
    \caption{Tonga Trench dataset: plain projection (left), 3D view (right).}
    \label{tonga_point}
\end{figure}

In our example we consider a dataset consisting of  8113 points.  We split it in a training and test set of 7000 and 1113 random samples without repetition (see Figure \ref{tonga_point}).\\

\begin{table}[h]
    \centering
    
\begin{tabular}{|*{7}{p{17mm}|}}

\cline{2-7}
\multicolumn{1}{c}{} & \multicolumn{3}{|c}{Gaussian kernel ($\varphi_1$)} &\multicolumn{3}{|c|}{Matérn kernel ($\varphi_2$)}\\
    \hline   
    $\tau$ &  time (s)  & RMAE  & RRMSE & time (s) & RMAE  & RRMSE \\
\hline
    1e-04 &   1.59e+03 &  6.99e-01 &  5.68e-02 &  1.34e+03 &  6.29e-01 &  5.45e-02 \\
\hline
   1e-05 & 1.57e+03 &  6.60e-01 &  6.14e-02 &  1.33e+03 &  6.16e-01 &  5.62e-02 \\
\hline

\end{tabular}
    \caption{Computational time,  RMAE and RRMSE using BO optimizer for Gaussian and Matérn kernels, $\varphi_1$ and $\varphi_2$,  using Tonga dataset. Two tolerances $\tau$ for the training error are used.}
    \label{tab:tonga}
\end{table}

\begin{figure}
    \centering
    \includegraphics[width=.49\textwidth]{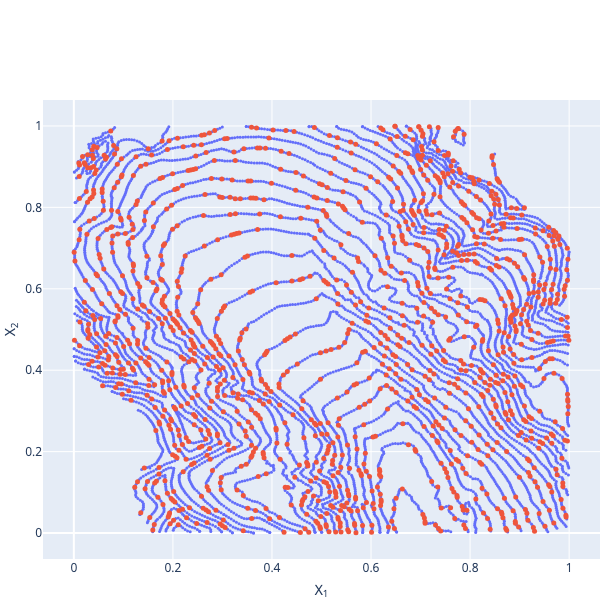} 
\includegraphics[width=.49\textwidth]{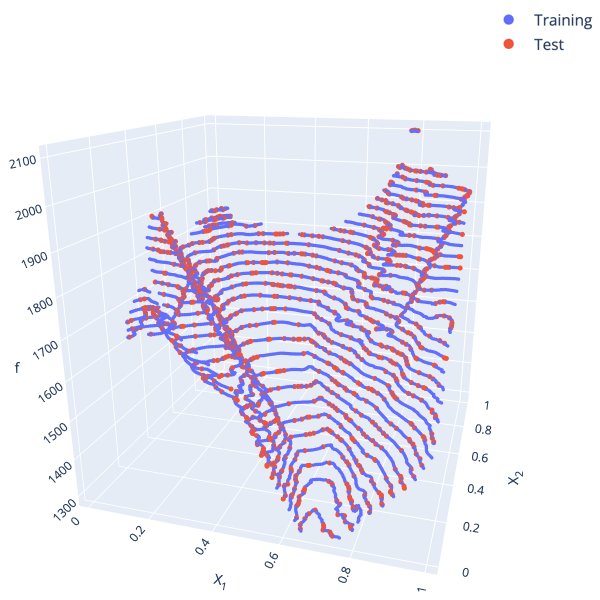}
    \caption{Franke's glacier dataset: plain projection (left), 3D view (right).}
    \label{glacier_points}
\end{figure}

\textbf{Franke's Glacier Dataset:} This dataset previously used in \cite{cav17} for interpolation of scattered data using RBFs for surface fitting, consists of 8338 measurements of altitude of a glacier.
Unfortunately we can not find any background on where these data were collected or indeed even the location of this glacier. More details on this dataset can be found at \url{https://search.r-project.org/CRAN/refmans/fields/html/glacier.html}. However, it is an interesting dataset in which it appears that the elevations are reported along lines of equal elevation, i.e. contours, perhaps from a digitization of a topographic map or survey.
In our example we consider the whole dataset consisting of 8338 points   and we split it in  a training and test set  of 7000 and 1338 random samples without repetition (see Figure \ref{glacier_points}).

\begin{table}
    \centering
    
\begin{tabular}{|*{7}{p{17mm}|}}

\cline{2-7}
\multicolumn{1}{c}{} & \multicolumn{3}{|c}{Gaussian kernel ($\varphi_1$)} &\multicolumn{3}{|c|}{Matérn kernel ($\varphi_2$)}\\
    \hline   
    $\tau$ &  time (s)  & RMAE  & RRMSE & time (s) & RMAE  & RRMSE \\
\hline
    1e-04 &   1.82e+03 &  9.25e-03 &  1.15e-03 &  1.62e+03 &  9.26e-03 &  8.74e-04 \\
\hline
   1e-05 & 1.83e+03 &  3.48e-02 &  1.49e-03 &  1.63e+03 &  9.33e-03 &  8.57e-04 \\
\hline

\end{tabular}
    \caption{Computational time,  RMAE and RRMSE using BO optimizer for Gaussian and Matérn kernels, $\varphi_1$ and $\varphi_2$,  using glacier dataset. Two tolerances $\tau$ for the training error are used.}
    \label{tab:glacier}
\end{table}


\subsubsection*{Acknowledgments}  

This research has been accomplished within the RITA \lq\lq Research ITalian network on Appro\-xi\-ma\-tion\rq\rq\ and the UMI Group TAA \lq\lq Approxi\-ma\-tion Theory and Applications\rq\rq. This work has been supported by the INdAM--GNCS 2022 Project \lq\lq Computational methods for kernel-based approximation and its applications\rq\rq, code CUP$\_$E55F22000270001, and by the Spoke ``Future HPC \& BigData'' of the ICSC--National Research Center in \lq\lq High-Perfor\-man\-ce Computing, Big Data and Quantum Computing", funded by European Union -- NextGenerationEU. Moreover, the work has been supported by the Fondazione CRT, project 2022 \lq\lq Modelli matematici e algoritmi predittivi di intelligenza artificiale per la mobilit$\grave{\text{a}}$ sostenibile\rq\rq.


%
%
%

\begin{thebibliography}{00}



\bibitem{acd} G. Allasia, R. Cavoretto, A. De Rossi, Hermite-Birkhoff interpolation on scattered data on the sphere and other manifolds. Appl. Math. Comput. 318, (2018) 35--50.

\bibitem{bab97} I. Babuška, J.M. Melenk, The partition of unity method, Internat. J. Numer. Methods Engrg. 40(4) (1997), 727--758.

\bibitem{bia18} M.E. Biancolini Fast Radial Basis Functions for Engineering Applications, Springer Cham, 2018.

\bibitem{Brochu} E. Brochu, V.M. Cora, N. De Freitas, A tutorial on Bayesian optimization of expensive cost functions, with application to active user modeling and hierarchical reinforcement learning, 2010, arXiv:1012.2599 

\bibitem{cav21a} R. Cavoretto, Adaptive radial basis function partition of unity interpolation: A bivariate algorithm for unstructured data, J. Sci. Comput. 87  (2021) 41.


\bibitem{Sandro} R. Cavoretto, A. De Rossi, S. Lancellotti, E. Perracchione, Software implementation of the partition of unity method, Dolomites Res. Notes Approx. 15 (2022) 35--46. 

\bibitem{cav21} R. Cavoretto, A. De Rossi, M.S. Mukhametzhanov, Ya.D. Sergeyev, On the search of the shape parameter in radial basis functions using univariate global optimization methods, J. Global Optim. 79 (2021) 305--327.

\bibitem{cav17}  R. Cavoretto, A. De Rossi, E. Perracchione, Optimal selection of local approximants in RBF-PU interpolation, J. Sci. Comput. 74 (2018) 1--22. 

\bibitem{cav22} R. Cavoretto, A. De Rossi, A. Sommariva, M. Vianello, RBFCUB: A numerical package for near-optimal meshless cubature on general polygons, Appl. Math. Lett. 125 (2022) 107704.

\bibitem{Fasshauer} G.E. Fasshauer, Meshfree Approximation Methods with MATLAB, World Scientific, Singapore, 2007.

\bibitem{Fasshauer15}
G.E. Fasshauer, M.J. McCourt, Kernel-based Approximation Methods Using MATLAB, World Scientific, Singapore, 2015.

\bibitem{Fornberg-Wright04} 
B. Fornberg, G. Wright, Stable computation of multiquadrics interpolants for all values of the shape parameter, Comput. Math. Appl. 47 (2004) 497--523.

\bibitem{fra18} E. Francomano, M. Paliaga, Highlighting numerical insights of an efficient SPH method, Appl. Math. Comput. 339 (2018) 899--915.


\bibitem{EI} D.R. Jones, M. Schonlau, W.J. Welch, Efficient global optimization of expensive black-box functions, J. Global Optim. 13 (1998) 455--492.


\bibitem{Larsson-Fornberg05} 
E. Larsson, B. Fornberg, Theoretical and computational aspects of multivariate interpolation with increasingly flat radial basis functions, Comput. Math. Appl. 49 (2005) 103--130.

\bibitem{laz02} 
D. Lazzaro, L. Montefusco, Radial basis functions for the multivariate interpolation of large scattered
data sets, J. Comput. Appl. Math. 140 (2002) 521--536. 

\bibitem{lin22} L. Ling, F. Marchetti, A stochastic extended Rippa’s algorithm for LpOCV, Appl. Math. Letters 129 (2022) 107955.

\bibitem{Lizotte} D. Lizotte, Practical Bayesian Optimization, PhD thesis, University of Alberta, Edmonton, Alberta, Canada, 2008.


\bibitem{Mockus_1978} J. Mockus, V. Tiesis, A. Zilinskas, The application of Bayesian methods for seeking the extremum, Towards Global Optimization 2 (1978) 117--129.

\bibitem{BOPy} F. Nogueira, Bayesian optimization: Open source constrained global optimization tool for Python, \url{https://github.com/fmfn/BayesianOptimization}

\bibitem{sklearn} F. Pedregosa, G. Varoquaux, A. Gramfort, V. Michel, B. Thirion, O. Grisel, M. Blondel, P. Prettenhofer, R. Weiss, V. Dubourg, J. Vanderplas,A. Passos, D. Cournapeau, M. Brucher, M. Perrot, E. Duchesnay, Scikit-learn: Machine learning in Python, Journal of Machine Learning Research 12 (2011) 2825--2830.

\bibitem{Rasmussen} C.E. Rasmussen, C. Williams, Gaussian Processes for Machine Learning, MIT Press, 2006.


\bibitem{shep68} D. Shepard, A two-dimensional interpolation function for irregularly-spaced data, in Proceedings of the 23rd National Conference ACM (1968) 517--523.

\bibitem{Practical} J. Snoek, H. Larochelle, R.P. Adams, Practical Bayesian optimization of machine learning algorithms, Advances in Neural Information Processing Systems 25 (2012) 2960--2968.

\bibitem{ren99} R. Renka, R. Brown, Algorithm 792: Accuracy tests of ACM algorithms for interpolation of scattered data in the plane, ACM Trans. Math. Softw. 25  (1999) 78--94.

\bibitem{scipy} P. Virtanen, R. Gommers, T.E. Oliphant, M. Haberland, T. Reddy, D. Cournapeau et al. SciPy 1.0: Fundamental Algorithms for Scientific Computing in Python, Nature Methods 17 (2020) 261--272.

\bibitem{wen02} H. Wendland, Fast evaluation of radial basis functions: methods based on partition of unity, in Approximation theory X: wavelets, splines and applications, Vanderbilt University Press, Nashville, (2002) 473--483.


\bibitem{wen05}
H. Wendland, Scattered Data Approximation, Cambridge Monogr. Appl. Comput. Math., vol. 17, Cambridge Univ. Press, Cambridge, 2005.


\end{thebibliography}
%

\end{document}